\documentclass[letterpaper, 10pt, conference]{ieeeconf}
\IEEEoverridecommandlockouts                              
\usepackage[utf8]{inputenc}
\usepackage{epsfig}

\usepackage{amsthm}
\usepackage{subfigure}
\usepackage{amsfonts}
\usepackage{amsmath}
\usepackage{eso-pic}
\usepackage{textcomp}
\usepackage{comment}
\usepackage{graphicx}
\usepackage{amssymb}
\usepackage{gensymb}
\usepackage[font=small,skip=0pt]{caption}
\usepackage[acronym]{glossaries}
\usepackage{xcolor}
\graphicspath{{Figures/}}
\usepackage{draftwatermark}
\makeatletter
\let\NAT@parse\undefined
\makeatother
\usepackage{hyperref}
\hypersetup{
colorlinks=true, %colorise les liens
breaklinks=true, %permet le retour à la ligne dans les liens trop longs
urlcolor= black, %couleur des hyperliens
citecolor=black,
linkcolor= blue, %couleur des liens internes
bookmarksopen=true,
pdftitle={Autonomous driving using an optimized neural network based adaptive LPV-MPC controller}, %informations apparaissant dans
pdfauthor={Y. Kebbati},
}

\title{\LARGE \bf
Coordinated PSO-PID based longitudinal control with LPV-MPC based lateral control for autonomous vehicles
}

\author{Yassine Kebbati$^{1*}$, Naima Ait-Oufroukh$^{1}$, Vincent Vigneron$^{1}$ and Dalil Ichalal$^{1}$% <-this % stops a space
\thanks{$^*$Corresponding author: \href{mailto:yassine.kebbati@univ-evry.fr}{yassine.kebbati@univ-evry.fr}}
\thanks{$^1$IBISC-EA4526, univ Evry, université Paris-Saclay, France}
}

\makeatletter
\newcommand*{\rom}[1]{\expandafter\@slowromancap\romannumeral #1@}
\makeatother

\begin{document}

\maketitle
%\conf{2022 European Control Conference (ECC)\\
%July 12-15, 2022. London, United Kingdom}
\thispagestyle{empty}
\pagestyle{empty}

%\tableofcontents

\SetWatermarkText{Preprint}

%%%%%%%%%%%%%%%%%%%%%%%%%%%%%%%%%%%%%%%%%%%%%%%%%%%%%%%%%%%%%%%%%%%%%%%%%%%%%%%%
\begin{abstract}
Autonomous driving is achieved by controlling the coupled nonlinear longitudinal and lateral vehicle dynamics. Longitudinal control greatly affects lateral dynamics and must preserve lateral stability conditions, while lateral controllers must take into account actuator limits and ride comfort. This work deals with the coordinated longitudinal and lateral control for autonomous driving. An improved particle swarm optimized PID (PSO-PID) is proposed to handle the task of speed tracking based on nonlinear longitudinal dynamics. An enhanced linear parameter varying model predictive controller (LPV-MPC) is also designed to control lateral dynamics, the latter is formulated with an adaptive LPV model in which the tire cornering stiffness coefficients are estimated by a recursive estimator. The proposed LPV-MPC is enhanced with an improved cost function to provide better performance and stability. Matlab/Carsim co-simulations are carried out  to validate the proposed controllers.

\end{abstract}

%%%%%%%%%%%%%%%%%%%%%%%%%%%%%%%%%%%%%%%%%%%%%%%%%%%%%%%%%%%%%%%%%%%%%%%%%%%%%%%%
\section{Introduction}

Autonomous vehicle technology has recently been the subject of significant research and development, due to its numerous advantages compared to conventional vehicles. Among its advantages are a reduced accident and fatality rate,  better and smoother traffic, energy and time efficiency. Automated driving systems consist of several modules that work together in a coordinated manner. First, there is a perception module made up  of a variety of sensors, including cameras, LIDAR, RADAR, GPS and IMU. Several algorithms are used by this module to detect the environment and extract relevant information. Second, a planning and decision-making module (speed and path planners) that exploits the perceived information to make decisions and plan future actions such as changes in  vehicle speed and  direction. Then comes the control system that executes the decisions taken by the previous module. Automatic control is the final step and one of the most important tasks in autonomous driving systems. The latter is divided into longitudinal and lateral control, where longitudinal control ensures precise speed tracking and lateral control handles the steering task. Extensive research has been done to address both tasks and reasonable progress has been made in the field. For instance, paper \cite{Marcano2018} developed an MPC controller for low velocity tracking in advanced driver-assistance systems (ADAS); it has been  compared to PID and neuro-fuzzy inference system (ANFIS) controllers. Simorgh \textit{et al.} \cite{simorgh2019adaptive} worked on adaptive PID control for speed regulation systems; they used the inverse model theory to adapt the PID gains. Good speed tracking has been obtained under adverse trajectory conditions and aerodynamic effects, but this approach depends on the reference model which is generally less accurate due to linearizations and simplified assumptions. Other works exploit the advantage of Artificial Intelligence to improve  classical  controllers such as Nie \textit{et al.}  \cite{nie2018longitudinal} who used radial basis function  to adapt the PID gains online (RBF-PID), which improved the tracking accuracy and driving comfort. In a similar way, Jin \textit{et al.} \cite{jin2020fuzzy} developed a fuzzy-PID for maintaining low speeds while driving downhill. Kebbati \textit{et al.}  \cite{kebbati2021optimized} developed a self-adaptive PID controller for the speed tracking task. They used genetic algorithms and neural networks to tune and optimize the PID gains for better tracking and disturbance rejection.

Several papers worked also on lateral control. In \cite{bujarbaruah2018adaptive}, an adaptive MPC controller has been designed for automatic lane-keeping. The developed controller was able to handle steering offsets learned from measured data by using a membership function approach. Overall, the results showed significant improvement provided that the vehicle lateral dynamics were perfectly known. Paper \cite{salt2021autonomous} also addressed lane-keeping, where an LPV-MPC was developed for steering control, and a dual rate extended Kalman filter (DREKF) was used for state estimation. The simulation results demonstrated acceptable performance. Likewise, Yang \textit{et al.} \cite{yang2020design} addressed the longitudinal control by designing a sliding mode controller (SMC) with a conditional integrator. They coupled it with an LPV-MPC for lateral guidance and then with an active disturbance rejection controller (ADRC). SMC suffers from the chattering phenomenon. ADRC showed good robustness against parametric uncertainties, but LPV-MPC achieved much more accurate tracking performance. Nevertheless, most of the previously mentioned papers only consider  constant longitudinal speeds, while variable speeds greatly affect lateral control performance. Study \cite{keb2021} developed an adaptive MPC for the path tracking task; they proposed an improved PSO algorithm to optimize the MPC parameters and used a lookup table approach to adapt these parameters online. Despite the optimal results for time-varying longitudinal velocities, this strategy cannot cover all possible situations and the lookup table method requires certain approximations which reduce the overall control precision. However, the same authors improved their controller design in \cite{kebb2021} by using neural networks and adaptive neuro-fuzzy inference systems to learn the optimal MPC parameters and replace the lookup table approach. The idea was to adapt the MPC controller to varying working conditions and external disturbances. The obtained results showed significant tracking improvements, but this approach still requires long offline optimizations. 

However, research on coupled longitudinal and lateral control is less abundant. Attia \textit{et al.} \cite{attia2014combined}  developed a combined lateral and longitudinal control strategy for autonomous driving. Nonlinear model predictive control (NLMPC) was used to manage lateral guidance and longitudinal control was based on Lyapunov theory. The control strategy achieved good results, but NLMPC requires very expensive computations. Although nowadays it can be applied in real-time, it requires expensive hardware. Chebly \textit{et al.} \cite{Chebly2019} also addressed  coupled lateral and longitudinal control. A multi-body formalism was used to accurately model the four-wheeled vehicle, and the control laws were developed based on Lyapunov function, then on the immersion and invariance with sliding mode approach. Robust speed and trajectory tracking were achieved with both methods. The immersion and invariance based controller showed better performance then its rival since it is less model-dependent. However, the design of this type of controllers is a difficult task and its practical implementation is not guaranteed. 

The contributions of this work are threefold; First, an improved LPV-MPC with a simple and precise prediction model is developed for handling the vehicle lateral dynamics. Second, recursive least squares algorithm (RLS) is used to estimate the cornering stiffness coefficients online and adapt the MPC prediction model iteratively. Third, an improved PSO algorithm is used to optimize a PID controller (PSO-PID) for handling the longitudinal dynamics which is coordinated with the lateral controller while ensuring lateral stability. This paper is organized as follows: Section \rom{2} presents the modeling of the longitudinal and lateral dynamics of the vehicle and the parameters estimation with recursive least squares algorithm. The design of the coordinated longitudinal and lateral control is detailed in section \rom{3}. Section \rom{4} presents and analyzes the obtained results. Finally, section \rom{5} concludes the work and gives some perspectives.  

\section{Vehicle Modeling}
The vehicle dynamics are highly nonlinear, where the longitudinal and lateral dynamics are coupled. In this paper a nonlinear longitudinal model that accounts for nonlinear power-train and tire longitudinal dynamics is developed and  coupled with a linear parameter varying bicycle model for the lateral dynamics.

\subsection{Longitudinal Dynamics Modeling} 
The longitudinal dynamics of a vehicle consist of several subsystems that represent the chassis, the engine, the transmission, the wheels, the tires and the brakes. Vehicle body dynamics can be modeled by the following equation \cite{kebbati2021optimized}:
\begin{equation}
\label{eq1}\small
F_t = m \frac{dv}{dt} + \frac{1}{2} \rho C_d A_f (v + v_w)^2 + m g C_r \cos \theta + m g \sin \theta
\end{equation}  
where  $F_t$, $m$, $v$ and $v_w$ are the drive force, the vehicle mass, the vehicle velocity and the wind velocity respectively. $A_f$, $C_d$ and $Cr$ are the vehicle cross sectional area, the drag and rolling resistance coefficients respectively.  $g$, $\rho$ and $\theta$ represent the  gravity, air density and road elevation angle. The propulsion system of electric vehicles is composed of the motor and the transmission. The motor can be modeled statically using the motor efficiency map and the transmission is modeled as a single speed gearbox as in equations (\ref{eq2}) and (\ref{eq3}) respectively \cite{kebbati2021optimized}:
\begin{equation}
\label{eq2}
\left\{
    \begin{array}{ll}
        T_e&= T_\text{e-ref}\\
        I_b&=\dfrac{T_e \omega_g}{u_b \eta^k}
    \end{array}
\right.
\end{equation}

\begin{equation}
\label{eq3}
\left\{
    \begin{array}{ll}
        T_g &= k_g T_e \\
        \omega_g &= k_g \omega_e
    \end{array}
\right.
\end{equation}
with $T_\text{e-ref}$, $T_e$ and $T_g$ being the engine reference torque (obtained from the efficiency map), the engine torque and the gear torque respectively. Similarly, $\omega_e$ and $\omega_g$ are the engine and the gear angular speeds. $k_g$, $I_b$ and $u_b$ represent the gear ratio, the current and the voltage drawn from the battery. $\eta$ is the motor efficiency where the exponent $k=1$ if $T_e \omega_g \geq 0$, otherwise $k=0$. The wheels can be modeled by equations (\ref{eq4}) under the assumption of negligible slip \cite{kebbati2021optimized}, and the tires are modeled by the magic formula given in (\ref{eq5}):
\begin{equation}
\label{eq4}
\left\{
    \begin{array}{ll}
        F &=\frac{T-B_d \omega_w +I_{\omega} \dot{\omega}_{\omega}}{R_{\omega}}\\
        \omega_w&= \frac{1}{R_w} v\\
        
    \end{array}
\right.
\end{equation}
\begin{equation}
\label{eq5}
F_i= F_z D \sin(C \tan^{-1}[{Bk-E[Bk-\tan^{-1}(Bk)]}])
\end{equation}
where $F$ and $T$ are either the drive or the braking force and torque respectively, $\omega_{\omega}$, $R_{\omega}$ and $v$ represent the wheel angular velocity, the wheel radius and the vehicle linear velocity. $I_{\omega}$ and $B_d$ are the wheel moment of inertia and damping coefficient. Equation (\ref{eq5}) is the Pacejka tire model \cite{Pacejka2008} where $F_i$ can be either the longitudinal or the lateral tire force, $F_z$ is the normal tire force, $k$ is the slip angle and $\{B,C,D,E\}$ are experimental fitting coefficients. The torque of the brake system is governed by equation (\ref{eq6}) where $B_a$ and $R_m$ are the actuator diameter and the brake pad mean radius. $P$ and $f_b$ are the brake pressure and the friction coefficient. These subsystems are modeled using the the powertrain-blockset of \texttt{MATLAB}.
\begin{equation}
\label{eq6}
T_b= \frac{f_b P \pi B_a^2 R_m}{2}
\end{equation}

\subsection{Lateral Dynamics Modeling}
The vehicle lateral dynamics are modeled using an LPV version of the standard bicycle model where the motion accounts for translation on the y-axis and rotation around the z-axis. By applying Euler-Newton's formalism to the model in (Fig. \ref{fig1}), using small angle approximations and linearized tire model, the full model can be summarized in the following equations \cite{keb2021}:

\begin{equation}
\label{eq7}
\left\{
    \begin{array}{ll}
        m(\dot{v}_y+v_x\Dot{\psi})&= 2[c_f(\delta-\frac{v_y + b\Dot{\psi}}{v_x}) + c_r\frac{a \Dot{\psi}- v_y}{v_x}]\\
        I_z \Ddot{\psi}&= 2[b c_f(\delta - \frac{v_y + b \Dot{\psi}}{v_x}) - a c_r \frac{b \Dot{\psi} - v_y}{v_x}]
    \end{array}
\right.
\end{equation}
where $v_{x/y}$, $\psi$ and $\delta$ are the longitudinal$/$lateral velocities, the heading angle and the steering angle respectively. $c_{f/r}$, $a$ and $b$ represent the cornering stiffness coefficients of front$/$rear tires and the distances from rear and front wheels to the vehicle gravity center (GC), and $I_z$ is the vehicle inertia moment. The corresponding affine LPV state space representation is given by: 
\vspace{1.5pt}
\begin{figure}[t]
\centering
\includegraphics[width=0.4\textwidth]{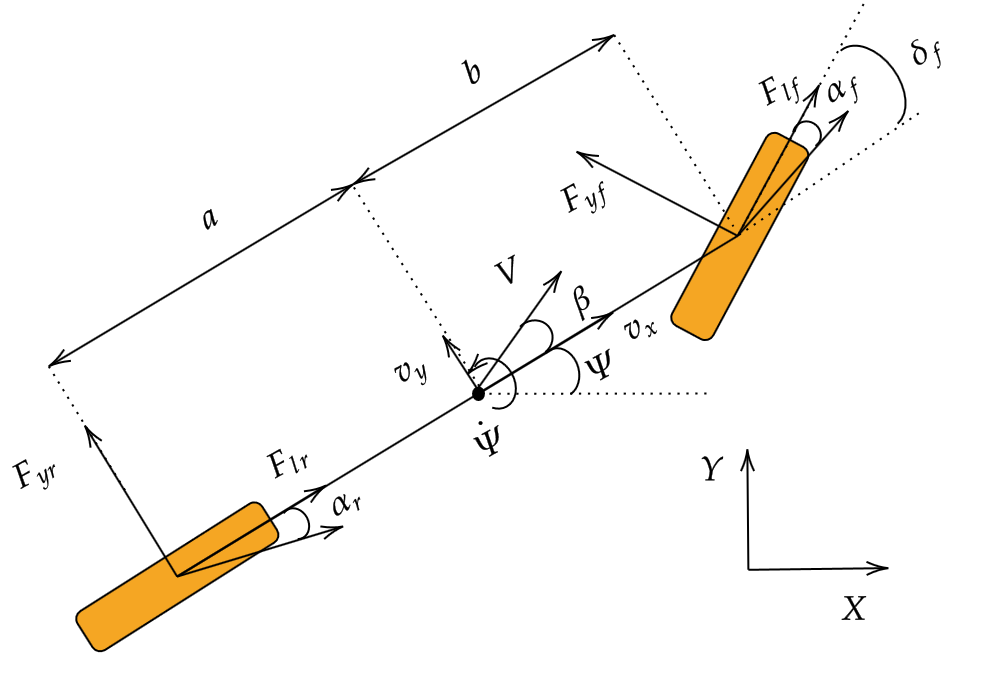}
\caption{Bicycle dynamic model}
\label{fig1}
\end{figure}

\begin{equation}
\label{eq8}
\left\{
    \begin{array}{ll}
        \Dot{X} &= A_d(\rho)X + B_d(\rho)u\\  %+Ed ==>$E = \left[\begin{array}{cccc} \frac{1}{4} & \frac{0.1}{10} & \frac{1.5}{10} & \frac{0.1}{50}\end{array}\right]^T,$
        Y &= CX
    \end{array}
\right.
\end{equation}
where $X=[\Dot{y}\ \psi\ \Dot{\psi}\ y]^T$, $Y=[y\  \psi]^T$ and $u=\delta_f$ represent the state vector, the output vector and the control signal respectively. The system matrices $A$ and $B$ depend on an LPV scheduling vector $\rho(t) = [v_x\ \frac{c_f}{v_x}\ \frac{c_r}{v_x}\  c_f]^T$. All the parameters are bounded and the system matrices are affine with respect to the parameters of the scheduling vector, the matrices are defined by:\\ \vspace{0.05cm}

$A = \left[{\begin{array}{cccc}
a_{11}(\rho)& 0 & a_{13}(\rho) & 0\\ 0 & 0 & 0 & 1\\a_{31}(\rho) & 0 & a_{33}(\rho) & 0\\1 & a_{42}(\rho) & 0 & 0\end{array}}\right],$\\\vspace{0.005cm}
$B = \left[\begin{array}{c}b_{11}(\rho) \\ 0 \\ b_{13}(\rho)\\ 0 \end{array}\right],$ \hspace{0.005cm}
$C = \left[\begin{array}{cccc} 0 &  0 & 0 & 1\\
0 & 1 & 0 & 0\end{array}\right]$.\\
\vspace{0.05cm}

\noindent where the terms of $A$ and $B$ are defined as follows:\vspace{0.1cm}

\noindent $a_{11}(\rho)=-2\frac{C_f + C_r}{m v_x}$, $a_{13}(\rho)=-v_x-\frac{bC_f-aC_r}{mv_x}$,\\ $a_{31}(\rho)=-2 \frac{bC_f-aC_r }{I_z v_x}$, $a_{33}(\rho)=-2\frac{C_f b^2+ C_r a^2}{I_z v_x}$,\\ $a_{42}(\rho)=v_x$, $b_{11}(\rho)=2\frac{C_f}{m}$, $b_{13}(\rho)=2\frac{C_fb }{I_z}$.\vspace{0.1cm}

\noindent Parameters $c_f$ and $c_r$ are obtained by using a recursive least squares (RLS) estimator. The estimator uses the linear tire cornering stiffness equations:

\begin{equation}
\label{eq9}
\left\{
    \begin{array}{ll}
        F_f &= c_f \alpha_f\\
        F_r &= c_r \alpha_r
    \end{array}
\right.
\end{equation}

The RLS algorithm given in equation (\ref{eq10}) uses the above mentioned equations as a parameter identification form where $z(k)$ is the measured value, $\phi(k)$ is a regression vector, $\theta(k)$ represents the estimated parameters and $e(k)$ is the identification error. These are defined in equations (\ref{eq11}): 

\begin{equation}
\label{eq10}
       z(k) =  \phi^T(k)\theta(k)+e(k)
\end{equation}

\begin{equation}
\label{eq11}
\left\{
    \begin{array}{ll}
        z(k) &= [F_f\ \ F_r]^T\\
        \phi^T(k) &= \left[\begin{array}{cc} \alpha_f &  0 \\
0 & \alpha_r \end{array}\right]\\
        \theta^T(k) &= [c_f\ \ c_r]^T
    \end{array}
\right.
\end{equation}

The values of $z(k)$ and $\phi(k)$ are obtained from the vehicle simulator (Carsim) and are fed to the RLS estimator, but in reality the lateral forces and side-slip angle are very hard to measure and they are usually estimated. The estimation process of these parameters is well studied in the literature \cite{lian2015cornering,baffet2009estimation} and is not in the scope of this study. Table ({\ref{tab1}}) summarizes the parameters used for both models.  

\begin{table}[htb]
\caption{Model parameters.} 
\label{tab1}
\centering
\begin{tabular}{c c c c c c} 
\hline
$m$ & $1575kg$ & $C_d$ & $0.29$ & $a$ & $1.6m$\\ [0.5ex] 

$C_r$ & $0.007$ & $k_g$ & $3.4$& $b$  & $1.2m$\\

$A_f$ & $1.6m^2$ & $R_w$ & $0.329m$ & $I_z$ & $2875N.m$\\

$\rho$ & $1.222kg/m^3$ & $R_m$ & $0.1778m$ & $I_{\omega}$  & $0.8kg.m^2$ \\

$f_b$  & $0.9$ & $B_d$ & $0.001/s^2$ & $\eta$ & $95\%$\\[1ex] 

\end{tabular}
\end{table}

\section{Controller Design}
The coupled lateral and longitudinal control strategy must preserve the lateral stability when performing lateral maneuvers. In this regard, both road geometry and lateral dynamics are to be considered \cite{attia2014combined}. The road information is used by the speed planner to estimate the maximum allowable longitudinal speed for the steering maneuver, the cruise speed has to be reduced in function of the curvature and must respect the following condition:
\begin{equation}
\label{eq12}
       v_\text{max}=\sqrt{\frac{g}{\rho_r}\frac{\phi_r+\mu}{1-\phi_r\mu}}
\end{equation}
with $g$, $\mu$, $\phi_r$ and $\rho_r$ being the gravity, the adhesion coefficient, the camber angle and the curvature respectively. This criterion is used as a predictive measure to guarantee admissible speeds before initiating lateral maneuvers. However, additional criteria related to lateral dynamics must be validated and in particular, the side slip angle ($\beta$) must obey criterion (\ref{eq13}):
\begin{equation}
\label{eq13}
       \beta \leq 10-\frac{7v_x^2}{40}
\end{equation}
The latter can be reformulated in terms of steering angle and implemented as a constraint on the control signal. According to \cite{Rajamani2012}, the side slip angle can be expressed by:

\begin{equation}
\label{eq14}
       \beta = \arctan \left(\frac{b\tan \delta_r + a\tan \delta_f}{a+b} \right)
\end{equation}
where the rear wheel steering angle can be set to zero as only the front wheels are steerable ($\delta_r=0$). Hence, using equation (\ref{eq14}), criterion (\ref{eq13}) can be reformulated into :

\begin{equation}
\label{eq15}
       |\delta_f| \leq \arctan \left(\frac{(a+b)tan(10-\frac{7v_x^2}{40})}{a} \right)
\end{equation}

\subsection{PSO-PID Control}

\begin{figure}[t]
\centering
\includegraphics[width=8.5cm,height=3.6cm]{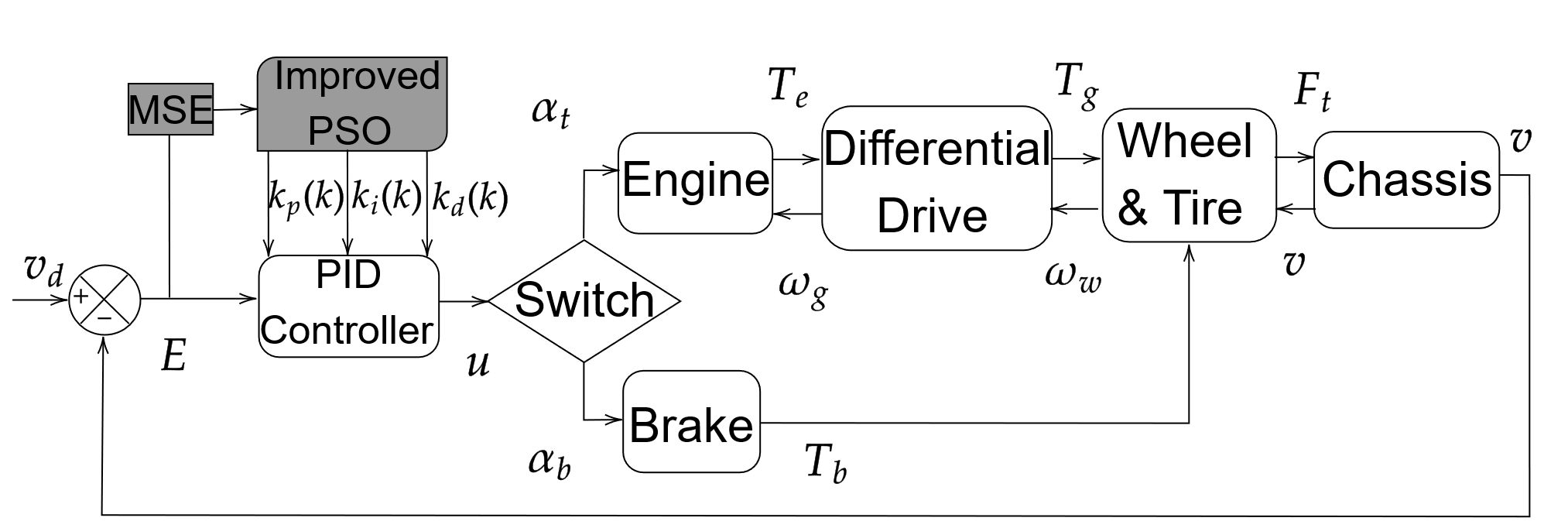}
\caption{PSO-PID approach.}
\label{fig2}
\end{figure}

For the speed tracking task, an optimized PID controller is proposed. PID control is commonly used in the industry, it achieves good performance whilst being relatively simple to design compared to other methods. The general formula of a PID controller is provided in equation (\ref{eq16}):
\begin{equation}
\label{eq16}
       u(t) = K_p e(t) + K_i \int e(\tau)d\tau + K_d \frac{de(t)}{dt}
\end{equation}
where $e(t)$ is the error between the reference and the output, and \{$K_p$, $K_i$, $K_d$\} are gains for the proportional, integral and derivative actions, which are often hard to tune. Therefore, the three parameters are optimized through an improved particle swarm optimization algorithm (PSO), which is a well established stochastic optimization technique \cite{song2021improved,xie2019novel}. The algorithm relies on the principle of swarm social behaviour, where a swarm contains multiple particles which represent potential solutions. Each particle has a position $p_i$ and a velocity $v_i$, and the idea is to find the velocity that moves particles towards the global optimal position. The algorithm is governed by:
\begin{equation}
\label{eq17}
\left\{
    \begin{array}{ll}
 v_i(k+1) =& \omega v_i(k) + c_1 r_1 (Pb_i(k)-x_i(k)) \\
 &+ c_2 r_2(Gb(k)-x_i(k))\\
 x_i(k+1)  =& x_i(k) + v_i(k+1)
    \end{array}
\right.
\end{equation}
where $\omega$, $c_1$ and $c_2$ are known as inertia weight, cognitive and social accelerations respectively, and $r_{1,2} \in [0,1]$ are random constants. $Gb$ and $Pb$ are the global best position of the whole swarm and the local best position in a single swarm generation. The inertia weight, social and cognitive accelerations are prefixed constants in the classic algorithm. However, in the improved version of this work they are dynamic according to equations (\ref{eq18},\ref{eq19}) which improves the PSO search performance: 
\begin{equation}
\label{eq18}
\omega = \omega_\text{min} + \frac{\exp{(\omega_\text{max}-\lambda_1(\omega_\text{max}+\omega_\text{min})\frac{g}{G})}}{\lambda_2}
\end{equation}
\begin{equation}
\label{eq19}
\left\{
    \begin{array}{ll}
 &c_1(k+1) = c_1(k) + \alpha \\
 &c_2(k+1) = c_2(k) + \beta \\
 &\alpha = -2\beta = 0.085\quad \text{for} \quad \frac gG \leq 30\%\\
 &\alpha = \frac{-\beta}{2} = 0.045\quad \text{for} \quad 30\% \leq \frac gG \leq 60\%\\
  &\alpha = \frac{-\beta}{2} = -0.025\quad \text{for} \quad 60\% \leq \frac gG \leq 85\%\\
 &\alpha = -\beta = -0.0025\quad \text{for} \quad \frac gG \geq 85\%
    \end{array}
\right.
\end{equation}
Parameters $g$ and $G$ are the actual and the last generations respectively. $\lambda_{1,2}$ are constants adjusted to ensure an exponential decrease from $\omega_\text{max}$ to $\omega_\text{min}$ which are maximum and minimum inertia weights. The advantage of this PSO improved version compared to other improved versions is that it enhances the overall search capabilities of the algorithm, the exponential decrease of $\omega$ accelerates the convergence towards the global best solution. On the other hand, increasing $c_1$ pulls the particles towards $Pb$ and enhances the exploration phase, while increasing $c_2$ speeds the convergence towards $Gb$ which enhances the exploitation phase and vice versa. The proposed PSO is used to optimize the PID gains \{$K_p$, $K_i$, $K_d$\} using the longitudinal model discussed in section \rom{2}. The approach is illustrated in (Fig. \ref{fig2}) where the mean squared error (MSE) is used as the fitness function of the PSO algorithm. The switch module in the figure decides whether the control signal is a braking or an acceleration command based on its sign, the switching logic sends the command either to the braking system or to the engine at a time \cite{nie2018longitudinal}. 

\subsection{LPV-MPC  Control}
\begin{figure}[t]
\centering
\includegraphics[width=8.5cm,height=5.5cm]{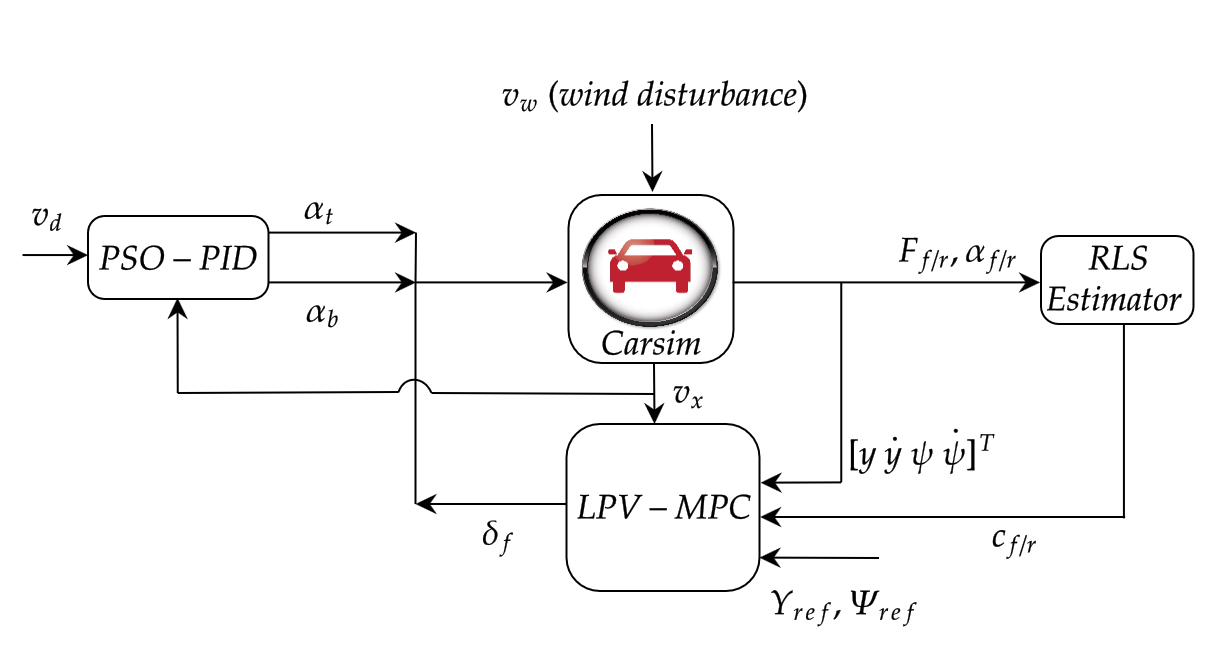}
\caption{Structure of the control approach}
\label{fig3}
\end{figure}

The lateral control performs the position and orientation tracking, an LPV-MPC formulation is proposed for this task. Model (\ref{eq8}) is discretized using sampling time $T_s$ and used as the MPC prediction model, where $A_d=I+T_sA$ and $B_d=T_sB$.  The principle relies on using the prediction model to predict the behaviour of the vehicle along a prediction horizon $N_p$, then solving the resulting optimization problem along a control horizon $N_c$. The goal is to find a control sequence that minimizes the error between the predicted vehicle behavior and the desired behaviour. At each iteration, a new adapted instance of the LPV model is used for the MPC predictions  which improves the overall prediction accuracy. Moreover, constraints are imposed on the control signal, its rate and system states to account for actuator limits and ride comfort. In addition to optimality, the possibility to handle constraints is a major advantage  of model predictive control. The standard optimization problem is given in equations (\ref{eq20}-\ref{eq24}):   

\begin{align}
\label{eq20}
min\ &J = (Y_r - Y_p)^T Q (Y_r - Y_p) + \Delta U^T R \Delta U\\
\label{eq21}
s.t:\ & x(k+1) = A_d(\rho) x(k) + B_d(\rho) \Delta u(k) \\
\label{eq22}
    & \Delta u_\text{min}\leq \Delta u(k) \leq \Delta u_\text{max}\\
\label{eq23}
    & u_\text{min}\leq u(k) \leq u_\text{max}\\
\label{eq24}
    & x_\text{min}\leq x(k) \leq x_\text{max}
\end{align}
where $\Delta U$, $Y_p$ and $Y_r$ represent the respective control sequence, the predicted and the reference trajectories, and $Q$ and $R$ are weighting matrices that penalize the tracking error and control effort. The constraints on the control and rate of control signal ensure a comfortable ride and impose actuator limits in addition to the lateral stability criterion given in equation (\ref{eq15}). To improve the optimization, an enhanced cost function is proposed by adding a relaxation factor ($\epsilon$) to avoid the infeasibility of strict constraints. An exponential weight ($\beta^{-k}$) is also included, the latter decreases with iteration $k$. This decreasing weight ensures that the first control input increments are more important in the sequence \cite{zhang2019electrical}. Therefore, immediate actions are favoured  to penalize immediate tracking errors more than future tracking errors since MPC is based on the receding horizon principle. The cost function is then reformulated into equation (\ref{eq25}):
\begin{equation}
\label{eq25}
\begin{array}{ll}
J = &\sum_{j=1}^{N_p}\beta^{-j}||y_r(k+j|k) - y_p(k+j|k)||^2_Q \\
&+\sum_{j=1}^{N_c}\beta^{-j}||\Delta u(k+j|k)||^2_R + \rho \epsilon^2
\end{array}
\end{equation}
The LPV-MPC is coded in \texttt{MATLAB} using \texttt{Yalmip} platform and solved using \texttt{Gurobi} optimization solver. The algorithm runs at $20Hz$ on an Asus Rog G17 with $2.6Ghz$ I7 $10750H$ and $32GB$ of DRAM. The coordinated control strategy is implemented and tested in \texttt{Matlab/Carsim} co-simulations as shown in Fig. \ref{fig3}, where $\alpha_t$ and $\alpha_b$ represent the throttle and brake signals, $v_d$, $Y_{ref}$ and $\Psi_{ref}$ are the desired longitudinal velocity, the lateral position and the heading trajectories. The desired velocity is generally obtained via a speed planner and must verify condition (\ref{eq12}). The LPV-MPC and PSO parameters are tuned iteratively until the desired performance is obtained. These are given in table (\ref{tab2}).  
\begin{table}[!h]
\caption{Model parameters.} 
\label{tab2}
\centering
\begin{tabular}{c c c c c c} 
\hline
$N_p$ & $9$ & $Q$ & $diag(35\ \ 3.25)$ & $R$ & $1.25$\\ [0.5ex] 

$\rho$ & $15$ & $\epsilon$ & $0.5$& $\beta$  & $3.5$\\

$T_s$ & $0.1s$ & $\Delta u_\text{max}$ & $\frac{pi}{12}$ & $\Delta u_\text{min}$ & $-\frac{pi}{12}$\\

$u_\text{max}$  & $\frac{pi}{6}$ & $u_\text{min}$ & $-\frac{pi}{6}$ & $G$ & $25$\\

$\omega_\text{max}$ & $1$ & $\omega_\text{min}$  & $0.1$ & $c_{i1}$ & $2.2$ \\

$\lambda_1$  & $3$ & $\lambda_2$ & $30$ & $c_{i2}$ & $2.2$\\\hline

\end{tabular}
\end{table}

\section{Results and Discussion}

\begin{figure}[t]
\centering
\includegraphics[width=8.5cm,height=5.5cm]{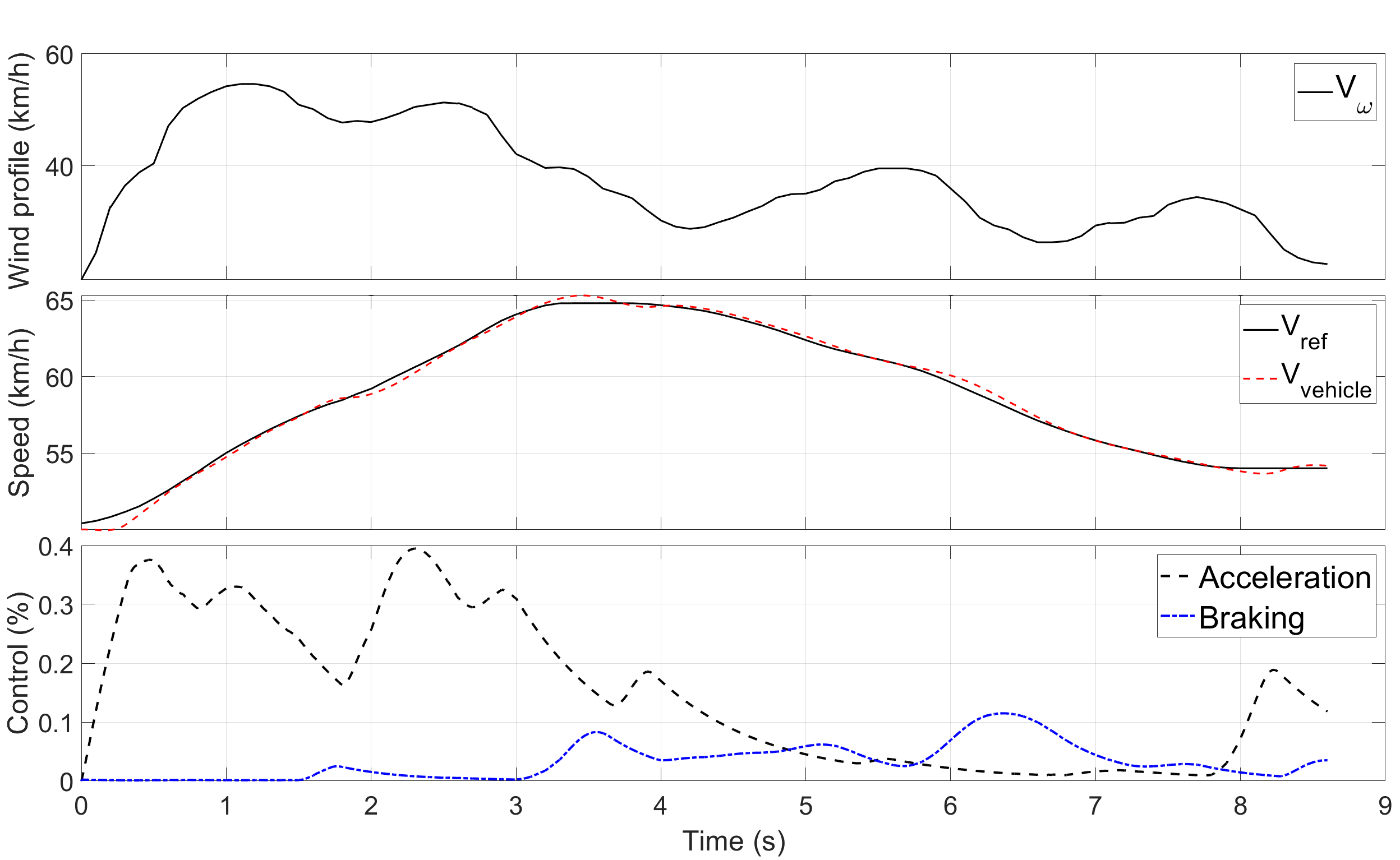}
\caption{Wind profile, speed tracking and control signals.}
\label{fig4}
\end{figure}

\begin{figure}[t]
\centering
\includegraphics[width=8.5cm,height=5.5cm]{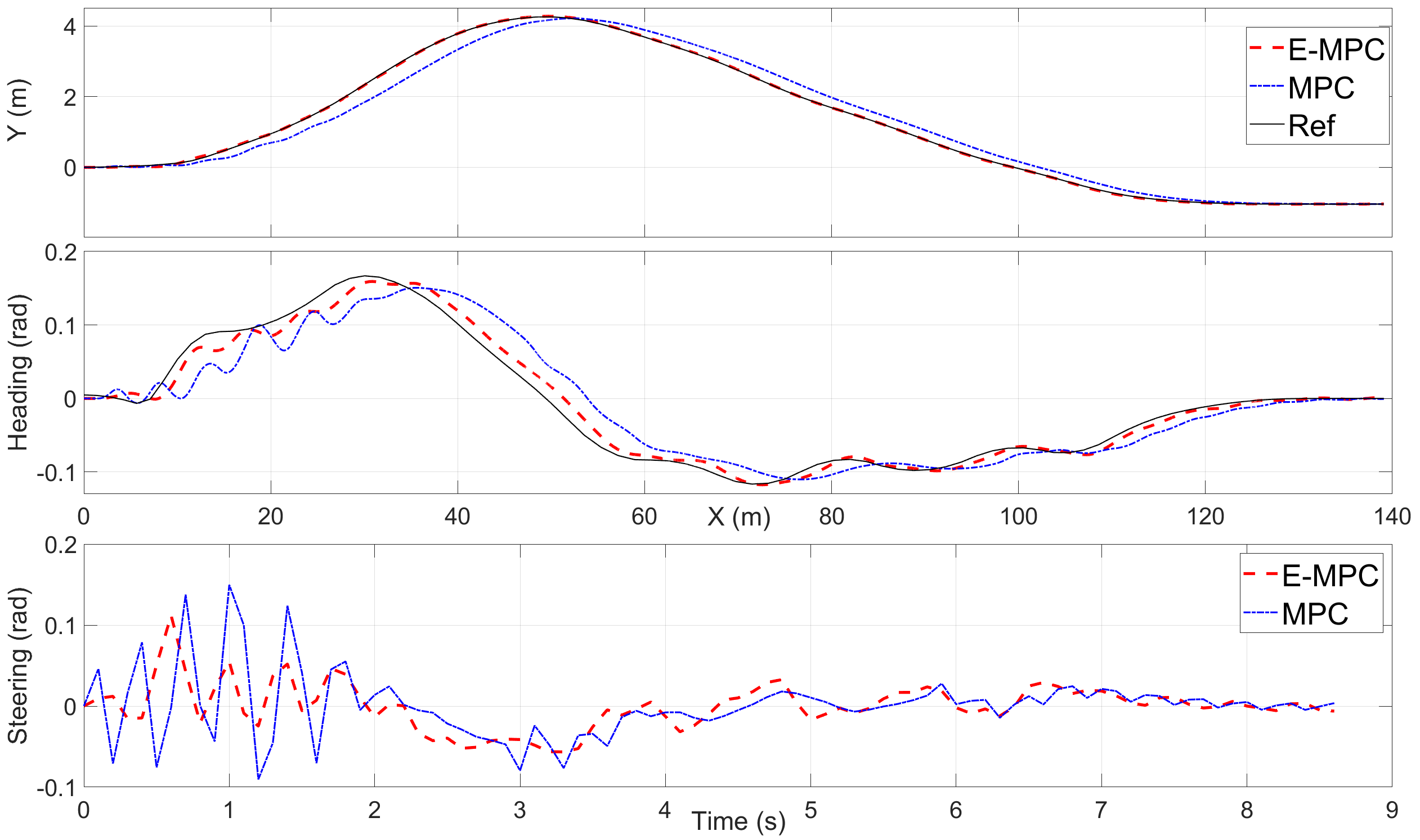}
\caption{LPV-MPC tracking performance and steering control.}
\label{fig5}
\end{figure}

The proposed control strategy is evaluated for a double lane-change maneuver where the velocity varies between $50$ and $65\ km/h$, and the vehicle is exposed to lateral and longitudinal wind disturbances. The test road is flat asphalt with and adhesion coefficient $\mu=0.95$. For this test, the PSO optimization generated the PID triplet \{$K_p=5.52$, $K_i=0.0547$, $K_d=0.469$\}. The LPV-MPC with regular cost function (\ref{eq20}) is compared to the one with the enhanced cost function (\ref{eq25}), they are denoted '\textit{MPC}' and '\textit{E-MPC}' respectively. Fig. \ref{fig4} shows the wind speed (with variable direction), the speed tracking performance and control signals. The results show smooth speed tracking with adequate accuracy, with an \textit{MSE=$0.0213$}. Fig. \ref{fig5} shows the lateral tracking performance, the advantage of \textit{E-MPC} is obvious in both position and yaw tracking with smoother steering which results in a more comfortable ride. The accuracy is measured by the \textit{MSE} value for position and orientation tracking, where \textit{E-MPC} scored $2.118e^{-4}$ against $0.601$ for \textit{MPC} with standard cost function. For the orientation tracking, \textit{E-MPC} scored $1.228e^{-4}$ against $6.756e^{-4}$. The highest position tracking error of \textit{E-MPC} does not exceed $5cm$. 

\begin{figure}[h]
\centering
\includegraphics[width=8.5cm,height=5cm]{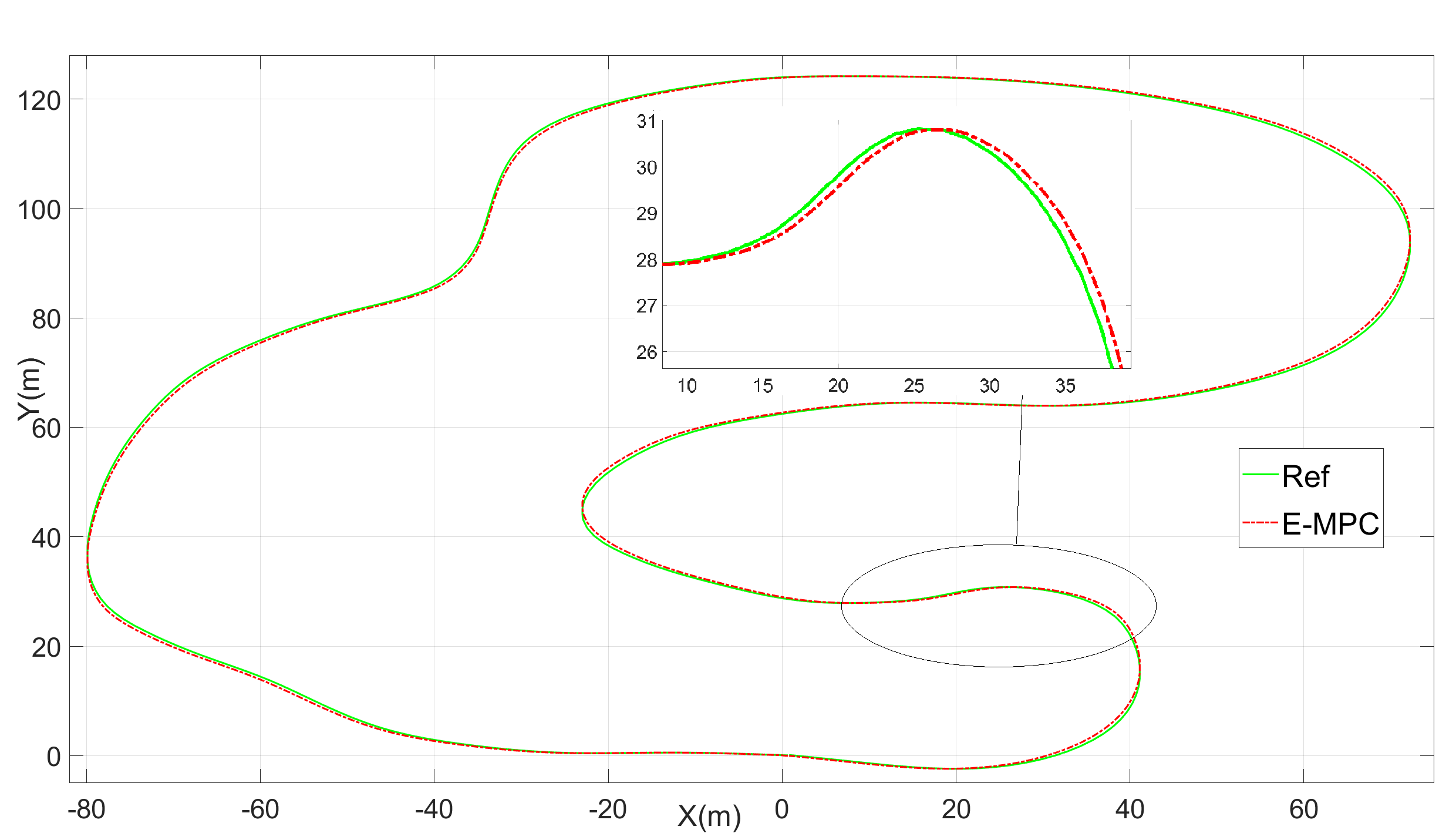}
\caption{Trajectory tracking.}
\label{fig6}
\end{figure}

\begin{figure}[h]
\centering
\includegraphics[width=8.5cm,height=5cm]{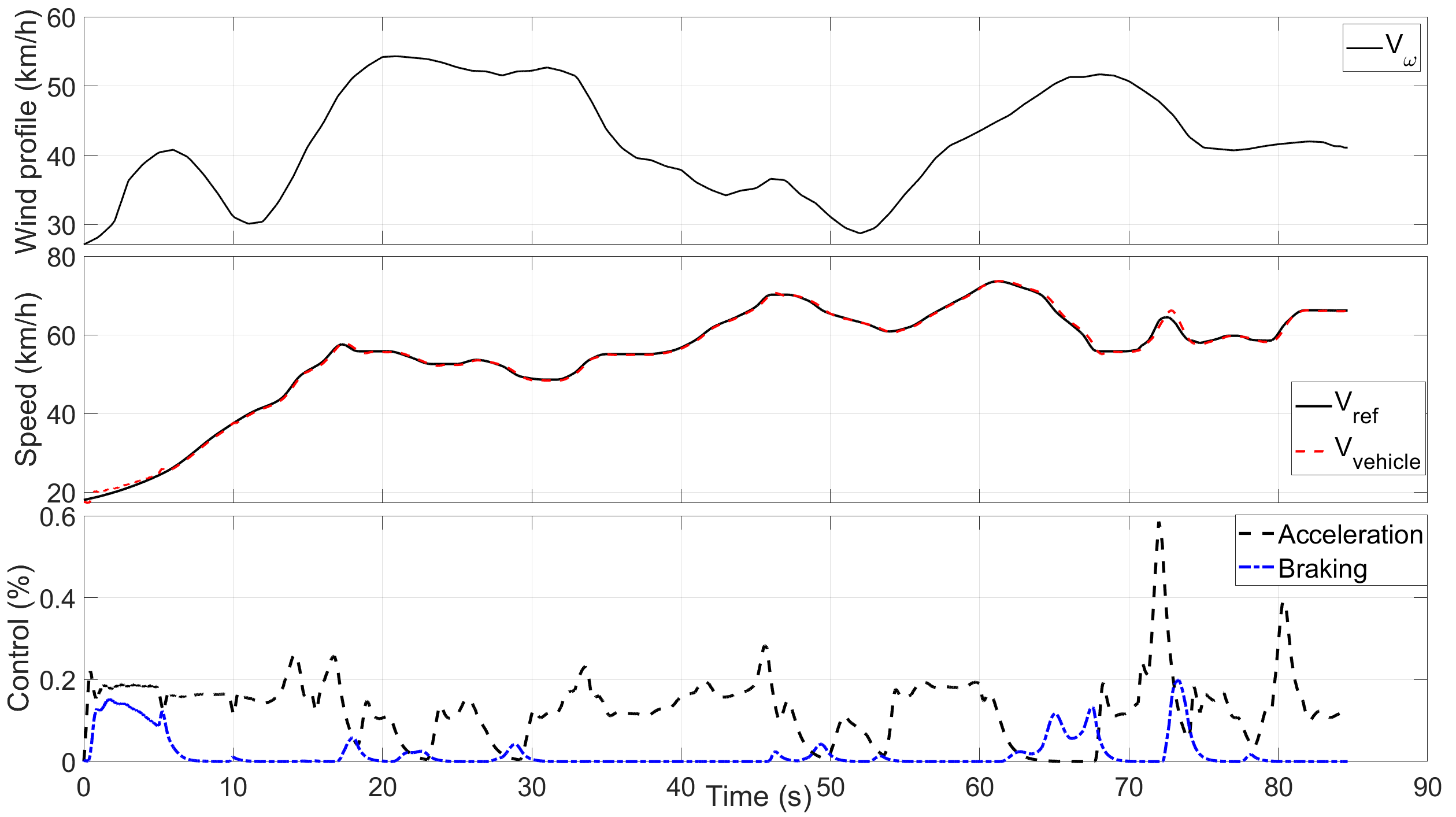}
\caption{Wind profile, Speed tracking and control signals.}
\label{fig7}
\end{figure}

\begin{figure}[h]
\centering
\includegraphics[width=8.5cm,height=5cm]{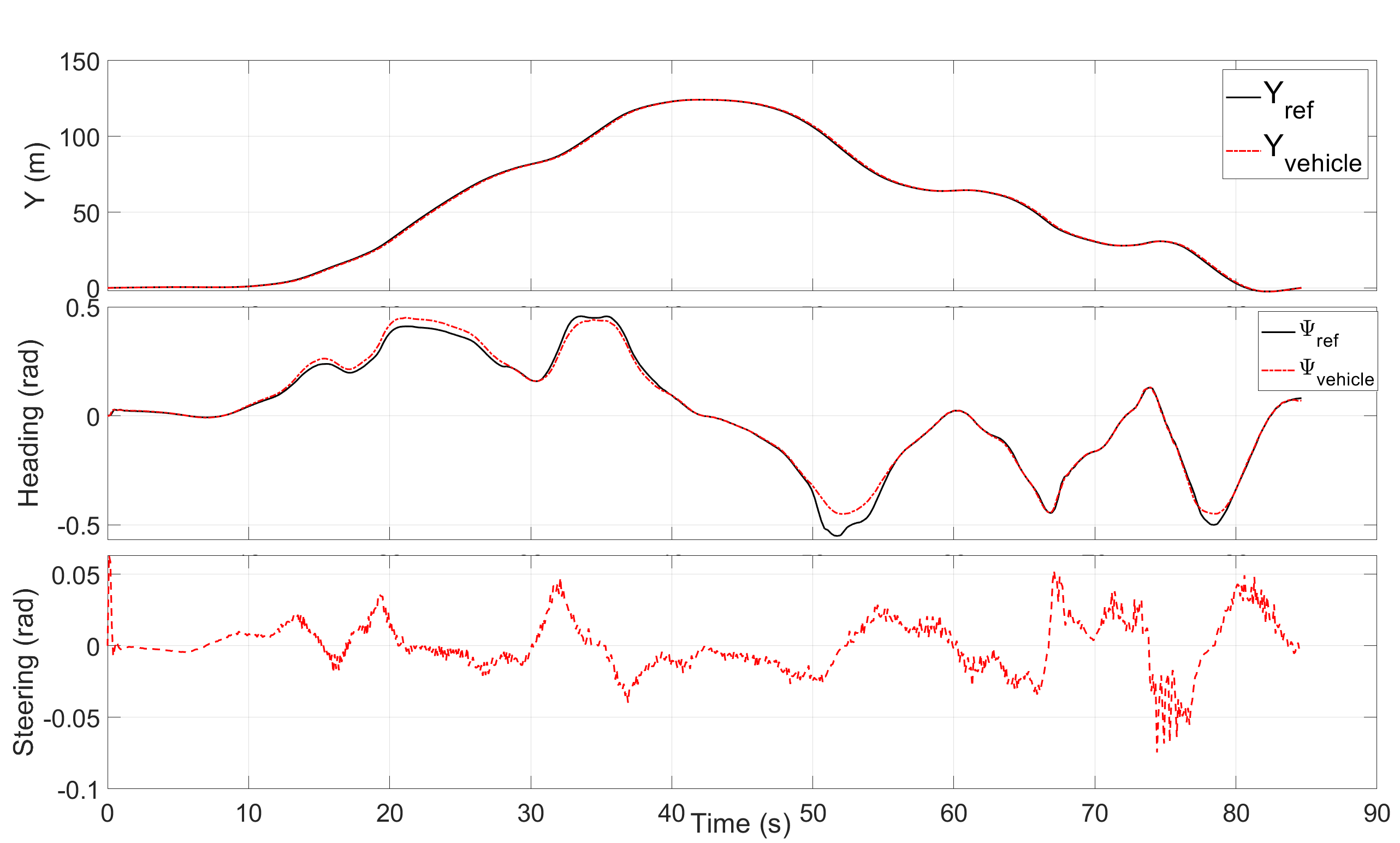}
\caption{Tracking performance and Steering control.}
\label{fig8}
\end{figure}

Furthermore, the effectiveness of \textit{E-MPC} is verified for a general trajectory (see Fig. \ref{fig6}) with different bends and a varying velocity profile. The same road characteristics were maintained and external wind disturbances were also applied. The PSO optimization in this test produced the PID triplet \{$K_p=1.884, K_i=0.048, K_d=0.325$\}. The wind profile, speed tracking performance and control signals are respectively shown in Fig. \ref{fig7}, where the observed performance is similar to the previous test, the tracking accuracy is indeed acceptable \textit{(MSE=$0.0187$)} and the controls are smooth. The respective lateral tracking performance along with the corresponding steering control are shown in Fig. \ref{fig8}. Despite the challenging trajectory and speed profile with imposed wind disturbances, the proposed LPV-MPC showed good performance overall and was able to track both position and orientation with a resulting \textit{MSE} of $0.32$ and $6.56e^{-4}$  for position and orientation tracking respectively. 
\section{Conclusions}
This paper addressed the coordinated longitudinal and lateral control for autonomous driving. A nonlinear model is used to account for power-train and vehicle longitudinal dynamics, and an optimized PSO-PID is developed to ensure accurate speed tracking, where the PID gains were optimized by a new improved PSO algorithm. An LPV-MPC with an enhanced cost function is designed for the lateral control task, the proposed controller takes into account the varying longitudinal velocity and cornering stiffness coefficients. A recursive least squares estimator is used to estimate the tire cornering stiffness coefficients. The proposed controllers were validated in \texttt{Matlab/Carsim} co-simulations for a double lane-change scenario and a general trajectory tracking problem. The results showed satisfactory performance for both lateral and longitudinal tracking while being relatively robust against wind disturbances. Future research shall address the development of an enhance LPV-MPC that handles both lateral and longitudinal control simultaneously.  

\bibliographystyle{ieeetr}
\bibliography{PhD}

\end{document}